# Twofold fast arithmetic

Evgeny Latkin

2014 Jul 8

**Abstract**: Can we assure math computations by automatic verifying floating-point accuracy? We define fast arithmetic (based on Dekker [1]) over twofold approximations $z \approx z_0 + z_1$, such that $z_0$ is standard result and $z_1$ assesses inaccuracy $\Delta z_0 = z - z_0$. We propose on-fly tracking $z_1$, detecting if $\Delta z_0$ appears too high. We believe permanent tracking is worth its cost. C++ test code for Intel AVX available via web.

Contents

History of this document

- 22-29 Jun 2014: complete reworking of earlier text by excluding stricter variants of arithmetic
- 8 Jul 2014: fixed a few misprints

## Motivation

Amazing progress of computers is only beginning yet. Mathematic computers should evolve to smarter, automatically identifying and addressing rounding problems, allow humans privilege be unaware of too much details. Let humans concentrate on their areas of interest, construction, science, education, etc.

Assuring mathematic software statically is best approach that works perfectly in rare special cases, like elementary functions, sine/cosine, etc. However generally, result's accuracy depends on specific input and is hard to predict statically. So we compute in hope: increase precision and pray this be enough.

But shall we try sort of "daemon" who would automatically check all floating-point operations with 2x higher precision, and signal if standard and 2x-precise results deviate too much? This cannot catch all problems; but if catches majority, this would anyway make us surer about results.

Our daemon should not overestimate deviations to avoid irrelevant panic. And we want it work in on-fly manner: assess deviation in parallel with main computations. Balancing cost versus reliability, we would like majority of programmers find our daemon affordable for typical applications.

We construct our "twofold daemon" basing on Dekker [1] technique of 1971 for 2x-precise arithmetic. We adapt it to modern Intel/AMD processors, and shift accents from increasing precision to estimating inaccuracy of original calculations. Following new factors enable Dekker arithmetic meet our goals:

- Unlike 1970[th], modern processors are much faster than memory. Typically, we can afford up to 10 extra operations while CPU is fetching data, with minimal damage for overall performance.
- Fast fused-multiply-add (FMA) operation critical for performance of Dekker arithmetic is widely available nowadays with Intel and AMD inexpensive processors targeted for mass market.

Twofold is approximation of a real value $z$ with unevaluated sum of floating-point numbers $z \approx z_0 + z_1$. Given $x = x_0 + x_1$ and $y = y_0 + y_1$ and $z = x \circ y$, twofold $z_0 + z_1$ keeps $z_0 = x_0 \circ y_0$ bitwise identical to standard result, and $z_1$ measures inaccuracy of $z_0$ by approximating $\Delta z_0 = z - z_0$.

Given chain of such calculations, resulting $z_1$ would assess rounding errors accumulated by $z_0$. If $z_1$ itself were accurate enough, $z_1$ must keep small comparing $z_0$ if the main 1x-precise computation converges. We want such test be very fast for checking in on-fly manner. Our performance criteria:

- 2x slower than standard, if large data fit in computer memory but not in CPU cache
- 10x slower than standard calculations, if small data pre-fetched in processor cache

Each twofold takes 2x memory, so slowdown cannot be less than 2x if data not in cache. Fitting under 2x implies very efficiently utilizing CPU idle time, doing all additional calculations while fetching data. About 10x criterion, note that modern implementations of quad type are typically 100x slower than double.

Twofold fast arithmetic assesses $\Delta z_0$ with minimum operations, and meets these 2x and 10x criteria.

Unlike intervals, twofolds do not tend overestimating $\Delta z_0$, so risk of irrelevant panic is low. Twofolds may miss problems: assess $\Delta z_0$ as fine while actually it is not. We accept this risk, as twofolds would catch majority of problems, so verifying would anyway assure results.

Below, Background section references the techniques we use. Algorithms defines arithmetic formulas. C++ test code explains experimental implementation. Performance discusses test results. Applications shows twofolds in work. Discussion summarizes our observations and proposes CPU improvements.

C++ experimental implementation available via web, free for academic and non-commercial use:

https://sites.google.com/site/yevgenylatkin/

# Background

- [Exact transforms]
- [Dekker arithmetic]

## Exact transforms

Hereon we base on floating-point "exact transforms" developed and used by many authors (see [1-6]). The transforms are decomposing a floating-point sum or product into result and tail, $a + b \rightarrow s + t$ or $a \times b \rightarrow p + e$. Here $s$ and $p$ are the correctly rounded floating-point results, $t$ and $r$ are exactly the rounding errors (unless $e$ suffers underflow and is additionally rounded or flushed to zero).

Transforming $a/b \rightarrow q, r$ such that $a = qb + r$ is also exact (unless the remainder $r$ underflows), and one can easily compute $r$ with fused-multiply-add (FMA) function available with modern processors. Correctly rounded $\mathrm{sqrt}(a)$ and FMA functions enable easily computing exact residual for square root, $\sqrt{a} \rightarrow c, d$ where $a = c^2 + d$ (unless $d$ underflows).

Let $a \oplus b$, $a \ominus b$, $a \otimes b$, $a \oslash b$, $\mathrm{sqrt}(a)$ be floating-point variants of basic arithmetic operations. We assume rounding correctly to nearest-even which is the usual standard mode. Let $\mathrm{fl}(a + b)$ be correctly rounded of $a + b$, and $\mathrm{err}(a + b)$ be its rounding error. With this notation we have $a \oplus b = \mathrm{fl}(a + b)$ and $\mathrm{err}(a + b) = a + b - a \oplus b$, and similarly for other operations.

Simple algorithms we use here proven for IEEE-754-2008 binary types and may be not valid for decimal or non-standard floating-point formats. Primarily we target the standard binary32 and binary64 formats, which correspond to the `float` and `double` types for majority of C/C++ implementations.

Let us take the following algorithms from Shewchuk [3]. Note that there is no if-then-else branching in these algorithms. This allows efficient vectoring for SIMD with modern processors, like Intel AVX.

Algorithm 1.1 (Fast Two-Sum): $a + b \rightarrow s + t$ provided $|a| \geq |b|$
(1) $s = a \oplus b$
(2) $b' = s \ominus a$
(3) $t = b \ominus b'$

Algorithm 1.2 (Two-Sum): $a + b \rightarrow s + t$ for arbitrary $a$ and $b$
(1) $s = a \oplus b$
(2) $b' = s \ominus a$
(3) $a' = s \ominus b'$
(4) $b^{\#} = b \ominus b'$
(5) $a^{\#} = a \ominus a'$
(6) $t = a^{\#} \oplus b^{\#}$

For subtracting $a - b \rightarrow d + t$ one could add $+(-b)$, but we prefer decomposing directly so saving one floating-point operation, Shewchuk [3]:

Algorithm 2.1: $a - b \rightarrow d + t$ provided $|a| \geq |b|$
(1) $d = a \ominus b$
(2) $b' = a \ominus d$
(3) $t = b' \ominus b$

Algorithm 2.2: $a - b \rightarrow d + t$ for arbitrary $a$ and $b$
(1) $d = a \ominus b$
(2) $b' = a \ominus d$
(3) $a' = b' \oplus d$
(4) $b^{\#} = b \ominus b'$
(5) $a^{\#} = a' \ominus a$

(6)  $t = a^\# \oplus b^\#$

We take the following obvious algorithm from Bailey and others [4]:

> <u>Algorithm 3 (Two-Product-FMA)</u>: $a \times b \to p + e$
(1)  $p = a \otimes b$
(2)  $e = \mathrm{fl}(ab - p)$

This algorithm uses the fused multiply-add (FMA) operation $\mathrm{fma}(a, b, -p) = \mathrm{fl}(ab - p)$. Supporting fast FMA with latest widely available processors makes other (slower and complicated) algorithms obsolete.

Floating-point dividing with exact remainder using FMA (see e.g.: Muller [6]):

> <u>Algorithm 4</u>: $a/b \to q, r$
(1)  $q = a \oslash b$
(2)  $r = \mathrm{fl}(a - qb)$

Floating-point square root with exact residual via FMA (Muller [6]):

> <u>Algorithm 5</u>: $\sqrt{a} \to c, d$
(1)  $c = \mathrm{sqrt}(a)$
(2)  $d = \mathrm{fl}(a - c^2)$

These floating-point exact transforms largely base on the following Sterbenz lemma, see Muller et al [6]:

> <u>Lemma 1 (Sterbenz)</u>: If floating-point $a$ and $b$ are close, so that $\frac{1}{2} \le a/b \le 2$, then their floating-point difference is exact, that is $a \ominus b = a - b$ exactly and $\mathrm{err}(a - b)$ is zero.

## Dekker arithmetic

Dekker [1] proposed simple and fast arithmetic on two-terms approximations like $z_0 + z_1$, which ideally can be up to 2x-precise. Dekker's arithmetic "renormalizes" results to ensure $z_0$ and $z_1$ do not overlap. We call renormalized pairs "coupled" numbers to distinguish from general-case twofolds.

Renormalizing means exact transform $z_0 + z_1 \to s + t$ with Algorithm 1.1 or 1.2 and replacing $z_0$ and $z_1$ with $s$ and $t$, so that renormalized $z_1$ would never exceed $\mathrm{ulp}(z_0)/2$ by magnitude. -- Hereon, $\mathrm{ulp}(u)$ is the "unit in last place" of a floating-point number $u$.

Dekker summation/subtraction and multiplication look like following. Assuming the basic floating-point operators $\oplus$ and $\ominus$ associate to left:

> <u>Algorithm 6</u>: $(x_0 + x_1) \oplus (y_0 + y_1) \to z_0 + z_1$
(1)  $z_0 = x_0 \oplus y_0$
(2)  $z_1 = x_1 \oplus y_1 \oplus \mathrm{err}(x_0 + y_0)$
(3)  Renormalize $z_0 + z_1$

> <u>Algorithm 7</u>: $(x_0 + x_1) \ominus (y_0 + y_1) \to z_0 + z_1$
(1)  $z_0 = x_0 \ominus y_0$
(2)  $z_1 = x_1 \ominus y_1 \oplus \mathrm{err}(x_0 - y_0)$
(3)  Renormalize $z_0 + z_1$

> <u>Algorithm 8</u>: $(x_0 + x_1) \otimes (y_0 + y_1) \to z_0 + z_1$
(1)  $z_0 = \mathrm{fl}(x_0 y_0)$
(2)  $z_1 = \mathrm{fl}(x_0 y_1) \oplus \mathrm{fl}(x_1 y_0) \oplus \mathrm{err}(x_0 y_0)$
(3)  Renormalize $z_0 + z_1$

With appropriate modifications, we reuse these algorithms for "fast" twofold arithmetic. Particularly we omit the renormalization step, which is irrelevant for our goals as it purges the useful information about inaccuracy of $z_0$. In turn, we cannot assume that input twofolds are non-overlapping.

We do not directly use Dekker's dividing and square root, instead propose better formulas leveraging fast FMA. Unlike early 1970[th] when original paper by Dekker [1] was published, very fast FMA is widely available nowadays with modern processors, so looks worth utilizing it.

This approach mainly targets extending C++ standard `double` type, as checking `float` by recalculating with `double` looks easier approach. In turn, twofold technique must work fine over quad-precision as basic type, so providing somewhat octal.

## Algorithms

- Methodic
- Add/subtract
- Multiply
- Divide
- Root
- 2x

### Methodic

We want twofolds assess inaccuracy accumulated by floating-point computations. Main part $z_0$ should remain bitwise same as for standard calculations, and $z_1$ should estimate accumulation of the rounding errors. If $z_1$ gets large comparing $z_0$, this should signal that precision of $z_0$ appears not enough. But let us allow misbalanced twofolds with $z_1$ large comparing $z_0$, let user's program react on such situations.

Given twofold $x = x_0 + x_1$ and $y = y_0 + y_1$ and operation $z = x \circ y$, define resulting $z_0 + z_1 \approx z$ like correctly rounded main part $z_0 = \mathrm{fl}(x_0 \circ y_0)$ and estimate $z_1 \approx \mathrm{fl}(\Delta z_0)$ for its deviation $\Delta z_0 = z - z_0$. Deviation of such estimate is $\Delta z_1 = \Delta z_0 - z_1 = z - (z_0 + z_1)$. Ideally, $z_1 = \mathrm{fl}(\Delta z_0)$ should be correctly rounded, so $z_0 + z_1$ is best possible twofold approximation, and $|\Delta z_1| \leq \mathrm{ulp}(z_1)/2$.

Unfortunately, such strict arithmetic would be slow. For monitoring $z_1$ in on-fly manner, we construct "fast" algorithms with minimal extra operations. For example, if we let $e_0 = \mathrm{err}(x_0 + y_0)$ for twofold summation, then $\Delta z_0 = x_1 + y_1 + e_0$. We "naïvely" let $z_1 = x_1 \oplus y_1 \oplus e_0$ to assess $\Delta z_0$ quickly. Such naïve arithmetic cannot guarantee $\Delta z_1$ small versus $\mathrm{ulp}(z_1)$, but looks good enough for our goal.

We define separate simplified faster algorithms for special case if $y_1 = 0$. Additionally we define faster algorithms for special case if input twofolds are non-overlapping "coupled" numbers. If necessary, one can renormalize results for "coupled" input so simulate 2x-precise arithmetic similar to Dekker [1].

Our algorithms do not require any special care for processing corner cases like NaN, infinity, and out of range. If such case happens, we rely on basic operations, which grant 0/0 result in NaN, etc.

Algorithms defined below in sub-sections Add/subtract, Multiply, Divide, and square Root. Enumerated like for example, Algorithm RTF1 "square Root Fast algorithm #1 for general-case Twofold arguments".

## Add/subtract

Given twofold $x = x_0 + x_1$ and $y = y_0 + y_1$, let us compute a reasonable twofold approximation like $z_0 + z_1 \approx z$ for the exact sum $z = x + y$. By definition, $z_0$ behaves exactly like the original 1x-precision summation. Following is "fast" algorithm, same as by Dekker [1] except we do not renormalize result. Except renormalization, this is literally the Algorithm 6 from [Background](#).

Here $x \oplus y$ would be floating-point summation:

> <u>Algorithm ATF1</u>: $(x_0 + x_1) \oplus (y_0 + y_1) \rightarrow z_0 + z_1$
> (1) $z_0 = x_0 \oplus y_0$
> (2) $z_1 = x_1 \oplus y_1 \oplus \mathrm{err}(x_0 + y_0)$

Important particular case is $y_1 = 0$, adding a single-length value $y_0$ to a twofold accumulator:

> <u>Algorithm ATF1.1</u>: $(x_0 + x_1) \oplus y_0 \rightarrow z_0 + z_1$
> (1) $z_0 = x_0 \oplus y_0$
> (2) $z_1 = x_1 \oplus \mathrm{err}(x_0 + y_0)$

And similarly for $z = x - y$. Here $x \ominus y$ is floating-point subtraction:

> <u>Algorithm STF1</u>: $(x_0 + x_1) \ominus (y_0 + y_1) \rightarrow z_0 + z_1$
> (1) $z_0 = x_0 \ominus y_0$
> (2) $z_1 = x_1 \ominus y_1 \oplus \mathrm{err}(x_0 - y_0)$

And subtraction, if $y_1 = 0$:

> <u>Algorithm STF1.1</u>: $(x_0 + x_1) \ominus y_0 \rightarrow z_0 + z_1$
> (1) $z_0 = x_0 \ominus y_0$
> (2) $z_1 = x_1 \oplus \mathrm{err}(x_0 - y_0)$

For general-case Algorithms ATF1 and STF1 we cannot guarantee if $\Delta z_1$ would be small comparing $z_1$. For example, if $\varepsilon = \frac{1}{2}\mathrm{ulp}(1)$, let $x_0 + x_1 = 1 - \varepsilon$ and $y_0 + y_1 = \varepsilon - \varepsilon^2$. For twofold fast summation, $z_0 = 1 \oplus \varepsilon = \mathrm{fl}(1 + \varepsilon) = 1$ and $z_1 = (-\varepsilon) \oplus (-\varepsilon^2) \oplus \varepsilon = (-\varepsilon) \oplus \varepsilon = 0$, while exact $\Delta z_0 = -\varepsilon^2$. Stricter algorithm might resolve this inaccuracy, but this is out of our scope for now.

Partial-case Algorithms ATF/STF 1.1 are strict and output the correctly rounded $z_1 = \mathrm{fl}(\Delta z_0)$.

There is no special algorithms for non-overlapping "coupled" input. The "fast" algorithms ATF/STF 1 and 1.1 are the best for this case. Because $z_0$ and $z_1$ are almost non-overlapping here, $z_0 + z_1$ provide nearly 2x-precise approximation of $z$. If you need non-overlapping output like Dekker summation, renormalize $z_0 + z_1$ with fast Algorithm 1.1. Fast renormalization requires 3 add/subtract operations.

Now let us assess performance of these algorithms by counting basic floating-point operations assuming exact transform $a + b \rightarrow s + t$ takes 6 operations:

|  | ATF/STF 1.1 | ATF/STF 1 |
|---|---|---|
| Exact transforms | 1 | 1 |
| More operations | 1 | 2 |
| Overall: | 7 | 8 |

Fast add/subtract algorithms meet the "slower by not more than 10x" performance criterion, as require less than 10 basic operations. They also meet the "2x slower if data not in CPU cache" criterion, see our testing results in the [Performance](#) section below.

## Multiply

Now let us compute a twofold approximation $z \approx z_0 + z_1$ for product $z = xy$ of twofold arguments $x = x_0 + x_1$ and $y = y_0 + y_1$. By definition $z_0 = x_0 \otimes y_0$ is exactly the single-length product of the main parts, so we have to compute $z_1$ for approximating $\Delta z_0 = z - z_0$.

Similarly to Dekker [1], our algorithms would base on the exact transforms $x_i \times y_j \to p_{ij} + e_{ij}$. If $e_{00}$ does not underflow, $\Delta z_0$ exactly equals $e_{00} + x_0 y_1 + x_1 y_0 + x_1 y_1$. Approximating $\Delta z_0$ naïvely, "fast" algorithm computes $z_1$ like floating-point sum of $e_{00}$ and corresponding $p_{ij}$.

Ordering of summation is not obvious however. To identify best formula, let us define "$\varepsilon$-order" of all involved terms, as measured by power of $\varepsilon$, where $\varepsilon = \frac{1}{2}\mathrm{ulp}(1)$ for the floating-point format we base. Let us classify relations of $x_1$ and $y_1$ versus $x_0$ and $y_0$ roughly like following:

(A) $x_0 + x_1$ and $y_0 + y_1$ are 2x-precise ("coupled"), so that $|x_1| \leq \varepsilon |x_0|$, and similarly for $y$
(B) $x_0 + x_1$ and $y_0 + y_1$ are more-or-less precise, like $|x_1| \leq \sqrt{\varepsilon}|x_0|$ or so, and similar for $y$
(C) $x_0$ and $x_1$ comparable so differ by less than by $\sqrt{\varepsilon}$ (e.g.: $x_0 = 1000$, $x_1 = -1$), and for $y$
(D) $x_0 + x_1$ form more-or-less precise inverse pair, so $|x_1| \geq 1/\sqrt{\varepsilon} |x_0|$, and similarly for $y$
(E) $x_0 + x_1$ form 2x-precise inverse "coupled", so that $|x_1| \geq 1/\varepsilon |x_0|$, and similarly for $y$

If we ignore possible cancellation of $e_{01} + e_{10}$ and $p_{01} + p_{10}$, following table summarizes the $\varepsilon$-orders comparing $z_0 = p_{00}$. For example, by design $|e_{00}| \leq \varepsilon |p_{00}|$, thus $\varepsilon$-order of $e_{00}$ comparing $p_{00}$ is $\varepsilon$:

|   | $e_{11}$ | $e_{01} + e_{10}$ | $e_{00}$ | $p_{11}$ | $p_{01} + p_{10}$ |
|---|---|---|---|---|---|
| A | $\varepsilon^3$ | $\varepsilon^2$ | $\varepsilon$ | $\varepsilon^2$ | $\varepsilon$ |
| B | $\varepsilon^2$ | $\varepsilon\sqrt{\varepsilon}$ | $\varepsilon$ | $\varepsilon$ | $\sqrt{\varepsilon}$ |
| C | $\varepsilon$ | $\varepsilon$ | $\varepsilon$ | 1 | 1 |
| D | 1 | $\sqrt{\varepsilon}$ | $\varepsilon$ | $1/\varepsilon$ | $1/\sqrt{\varepsilon}$ |
| E | $1/\varepsilon$ | 1 | $\varepsilon$ | $1/\varepsilon^2$ | $1/\varepsilon$ |

According to this table, the terms $e_{01}, e_{10}, e_{11}$ are anyway minor. If we ignore them, the table hints ordering from lower to higher magnitude summands like following:

(A) $e_{00}$ $+$ $(p_{01} + p_{10})$
(B) $e_{00} + p_{11} + (p_{01} + p_{10})$
(C) $p_{11} + (p_{01} + p_{10})$
(D) $p_{11} + (p_{01} + p_{10})$
(E) $p_{11}$

Formula $z_1 = e_{00} \oplus p_{11} \oplus (p_{01} \oplus p_{10})$ covers all these cases without too much of extra computations. Basing on this formula, "fast" algorithm for general-case twofolds would look as follows. In the step (1), we implicitly omit those computations that do not contribute to the result:

    <u>Algorithm MTF1</u>: $(x_0 + x_1) \otimes (y_0 + y_1) \to z_0 + z_1$
(1) $x_i \times y_j \to p_{ij} + e_{ij}$
(2) $z_0 = p_{00}$
(3) $z_1 = e_{00} \oplus p_{11} \oplus (p_{01} \oplus p_{10})$

For important partial case $y_1 = 0$, we can omit $p_{11}$ and $p_{01}$ which are zero:

    <u>Algorithm MTF1.1</u>: $(x_0 + x_1) \otimes y_0 \to z_0 + z_1$
(1) $x_i \times y_j \to p_{ij} + e_{ij}$
(2) $z_0 = p_{00}$
(3) $z_1 = e_{00} \oplus p_{10}$

Special case if input $x$ and $y$ are non-overlapping "coupled" twofolds allows omitting $p_{11}$, which in this case is minor comparing other terms of the summation. Following is special fast algorithm for "coupled" inputs. This algorithm is same as Dekker's except we omit renormalizing the result:

   <u>Algorithm MPF1</u>: $(x_0 + x_1) \otimes (y_0 + y_1) \to z_0 + z_1$
   (1) $x_i \times y_j \to p_{ij} + e_{ij}$
   (2) $z_0 = p_{00}$
   (3) $z_1 = e_{00} \oplus (p_{01} \oplus p_{10})$

For partial case if $y_1 = 0$, "coupled" algorithm MPF 1.1 would be literally same as MTF 1.1

Now let us assess performance of these algorithms by counting required basic operations. We count add/subtract, multiply, and FMA operations separately. Recall that each exact transform $a + b \to s + t$ takes 6 add/subtract operations, $a \times b \to p + e$ takes 1 multiply plus 1 of FMA (plus maybe 1 negation which we ignore). We count all operations, explicit and hidden in nested $n$-fold summations.

|  | Twofold | | Coupled | |
|---|---|---|---|---|
|  | Full | 1.1 | Full | 1.1 |
| $a + b$ | 3 | 1 | 2 | 1 |
| $a \times b$ | 3 | 1 | 2 | 1 |
| $a \times b \to p + e$ | 1 | 1 | 1 | 1 |
| Overall: | 8 | 4 | 6 | 4 |
| Add/sub: | 3 | 1 | 2 | 1 |
| Multiply: | 4 | 2 | 3 | 2 |
| FMA: | 1 | 1 | 1 | 1 |

Fast algorithms ATF/APF 1 and 1.1 must meet the "10x" criterion as takes less than 10 basic operations. See also <u>Performance</u> section below.

## Divide

Consider dividing twofold numbers $x = x_0 + x_1$ and $y = y_0 + y_1$. We a reasonably good approximation $z \approx z_0 + z_1$ for quotient $z = x/y$. By definition, $z_0$ equals correctly rounded $\text{fl}(x_0/y_0)$, and $z_1$ should approximate deviation $\Delta z_0 = z - z_0$.

Let us start with standard single-length division. Given floating-point $a$ and $b$ such that $b \neq 0$, consider iterative process for decomposing $a/b$ into: partial quotient $Q_N = q_0 + \cdots + q_n$ and remainder $r_{N+1}$ that approximates $a - Q_N b$. Here we leverage of fast fused-multiply-add (FMA) operation:

> <u>Process 1</u>: Decompose $a/b$ into $Q_N = q_0 + \cdots + q_N$ and $r_{N+1}$
> (1)  $r_0 = a$
> (2)  $q_n = \text{fl}(r_n/b)$
> (3)  $r_{n+1} = \text{fma}(r_n - q_n b)$

Steps (2) and (3) is the dividing with exact remainder Algorithm 4 from <u>Background</u>. Provided $r_{n+1}$ is not additionally rounded due to underflow of $q_n b$, well-known fact is that $r_{n+1}$ equals $r_n - q_n b$ exactly, see Muller [6]. So, if we ignore possible underflow, $r_{n+1}$ is exact remainder and $Q_N$ converges to exact $a/b$.

Additionally please note, that $q_n$ do not overlap each other, namely $|q_{n+1}| \leq \text{ulp}(q_n)/2$. Indeed, as $r_n$ is remainder of $r_{n-1}/b$, then $|r_n/b| \leq \text{ulp}(q_n)/2$ as otherwise dividing $r_{n-1}/b$ would result differently. Rounding $r_n/b$ cannot make magnitude of $q_{n+1} = \text{fl}(r_n/b)$ higher than $\text{ulp}(q_n)/2$, provided $\text{ulp}(q_n)/2$ is representable as floating-point number.

Therefore, $Q_1 = q_0 + q_1$ is best 2x-precise approximation for $a/b$ with correctly rounded $q_0$ and $q_1$.

This gives us following algorithm of twofold dividing $x_0 + x_1$ by $y_0 + y_1$ in special case if $x_1 = y_1 = 0$. Here we do not actually compute $r_2$ which do not contribute to result:

> <u>Algorithm DTF 1.1.1</u>: $x_0 \oslash y_0 \to z_0 + z_1$
> (1)  $z_0 = \text{fl}(x_0/y_0)$
> (2)  $r_1 = \text{fma}(x_0 - z_0 y_0)$
> (3)  $z_1 = \text{fl}(r_1/y_0)$

Note how this algorithm processes corner cases. If occasionally $y_0 = 0$, then automatically $z_0$ is infinite or NaN depending on $x_0$, thus so is $z_1$. If $x_0$ is occasionally NaN, this NaN propagates to z, etc. Thus, we do not need explicitly processing special arguments.

Now consider dividing a twofold $a = a_0 + a_1$ by 1x-precision $b$. Again, consider iterative decomposing into partial quotient and exact remainder. In this case, remainder would be twofold $r_n = r_{n0} + r_{n1}$:

> <u>Process 2</u>: Decompose $(a_0 + a_1)/b$ into $Q_N = q_0 + \cdots + q_N$ and $r_{N+1}$
> (1)  $r_0 = a$
> (2)  $q_n = \text{fl}(r_{n0}/b)$
> (3)  $r_{n+1,0} = \text{fma}(r_{n0} - q_n b)$
> (4)  $r_{n+1,1} = r_{n,1}$
> (5)  Renormalize $r_{n+1}$

Similar to Process 1 above, here $r_{n+1}$ is exact remainder, and $Q_N$ converges to exact $a/b$.

We cannot claim if $q_0 + q_1$ is non-overlapping 2x-precise. However, $q_2$ is anyway small comparing $q_1$, so $q_0 + q_1$ is still good approximation. For example, if $a_0 = a_1 = 1$ and $b = 1$, then $q_0 = q_1 = 1$ with exact reminder $r_2 = 0$.

This implies following algorithm of twofold dividing in special case if $y_1 = 0$, while $x_1$ may be non-zero. Here we explicitly unroll the Process 2 for this specific case and omit needless $r_2$ and $r_{11}$. Note, that this algorithm degrades to DFT 1.1.1 in case if $x_1 = 0$:

<u>Algorithm DTF 1.1</u>: $(x_0 + x_1) \oslash y_0 \to z_0 + z_1$
(1) $z_0 = \text{fl}(x_0/y_0)$
(2) $r = \text{fma}(x_0 - z_0 y_0)$
(3) $c = r \oplus x_1$
(4) $z_1 = \text{fl}(c/y_0)$

Finally, let us consider dividing twofolds $a = a_0 + a_1$ and $b = b_0 + b_1$. We consider two steps defining $q_0 + q_1$ and remainder. Twofold remainder is inexact, so we cannot recommend this process for $q_2$ etc.

- Let $q_0 = \text{fl}(a_0/b_0)$
- Let $r_0 = a_0 - q_0 b_0$ (note: $r_0$ is exact)
- Let $c_0 + c_1 + c_2 + c_3 = r_0 + a_1 - q_0 b_1$
- Let $d_0 + d_1$ be renormalized $b_0 + b_1$
- Let $q_1 = \text{fl}(c_0/d_0)$

With $c_0 + c_1 + c_2 + c_3$ we refer to known expand-and-distill technique used by many authors for multi-precision calculations, see [2-6]. With this techniques we could expand $q_0 b_1 \to p + e$ and compute sum of $r_0 + a_1 - (p + e)$ exactly resulting in $n$-fold $c = c_0 + \cdots$ with not more than four non-overlapping $c_i$. Additionally we assume $c_0 \neq 0$, if $c$ is not zero.

By design, $c$ is exact remainder of $q_0$, so that $a = q_0 b + c$ exactly. Thus best for $q_1$ would be rounded of $q = c/d$ if we could compute it. Instead we approximate with $q_1 \approx c_0/d_0$, quotient of the main parts of $c$ and $d$. Such approximation is accurate modulo approximately $3\varepsilon|q_1|$ where $\varepsilon = \text{ulp}(1)/2$. Indeed:

Suppose $d = d_0(1 + \delta)$ and $c = c_0(1 + \gamma)$ with some $\delta$ and $\gamma$ smaller than $\varepsilon = \text{ulp}(1)/2$. Then $c/d = c_0/d_0 \cdot (1 + \gamma)/(1 + \delta)$, so deviation $q - c_0/d_0 = c_0/d_0 \cdot ((1 + \gamma)/(1 + \delta) - 1)$. This coefficient magnitude $(1 + \gamma)/(1 + \delta) - 1 = (\gamma - \delta)/(1 + \delta)$ does not exceed $2\varepsilon/(1 - \varepsilon)$ or approximately $2\varepsilon$. Rounding $q_1 = \text{fl}(c_0/d_0)$ may add up to $\varepsilon|c_0/d_0|$ of inaccuracy. Thus overall distance $|q - q_1|$ does not exceed $3\varepsilon|c_0/d_0| \approx 3\varepsilon|q_1|$.

For our "fast" algorithm, we use this $q_1$, but simplify calculations like follows. If remainder $r_0 + r_1$ were exact, $c_0 + c_1$ would be best non-overlapping approximation for $c$. In practice, this is very good and fast formula for $c_0$ if we skip computing $c_1$ which we actually do not need:

- $r_0 = \text{fma}(a_0 - q_0 b_0)$
- $r_1 = \text{fma}(a_1 - q_0 b_1)$
- $r_0 + r_1 \to c_0 + c_1$

In overall, we come to the following "fast" algorithm for twofold dividing. Note that algorithm DTF 1.1 is degenerate variant of this algorithm for the case if $y_1 = 0$:

<u>Algorithm DTF 1</u>: $(x_0 + x_1) \oslash (y_0 + y_1) \to z_0 + z_1$
(1) $z_0 = \text{fl}(x_0/y_0)$
(2) $r_0 = \text{fma}(x_0 - z_0 y_0)$
(3) $r_1 = \text{fma}(x_1 - z_0 y_1)$
(4) $c_0 = r_0 \oplus r_1$
(5) $d_0 = y_0 \oplus y_1$
(6) $z_1 = \text{fl}(c_0/d_0)$

Simplified special case if input is non-overlapping "coupled". Here we can omit renormalizing $y_0 + y_1$:

<u>Algorithm DPF 1</u>: $(x_0 + x_1) \oslash (y_0 + y_1) \to z_0 + z_1$
(1) $z_0 = \text{fl}(x_0/y_0)$
(2) $r_0 = \text{fma}(x_0 - z_0 y_0)$
(3) $r_1 = \text{fma}(x_1 - z_0 y_1)$

(4)  $c_0 = r_0 \oplus r_1$

(5)  $z_1 = \text{fl}(c_0/y_0)$

"Coupled" variants for cases if $y_1 = 0$ and if $x_1 = y_1 = 0$ would literally repeat DTF 1.1 and 1.1.1

Because DPF1 saves just one basic summation of $y_0 \oplus y_1$, this must not make it much faster than DTF1. However, the point of DPF1 is that output $z_0 + z_1$ is easy to renormalize for simulating Dekker dividing. Here resulting $|z_1|$ is small comparing $|z_0|$ so fast renormalization algorithm 1.1 works fine.

Let us demonstrate how twofold dividing works:

<u>Example</u>: Given a floating-point format, let $\varepsilon = \text{ulp}(1)/2$. Then let $x_0 = 1$ and $x_1 = 0$, and $y_0 = 1 - \varepsilon$ and $y_1 = \varepsilon$. This way, $x/y$ is dividing 1/1, though with slight inaccuracy in $y_0$. With these data, we have $z_0 = \text{fl}(1/(1-\varepsilon)) = \text{fl}(1 + \varepsilon + \varepsilon^2 + \cdots) = 1 + 2\varepsilon$, thus $r_0 = \text{fma}\big(1 - (1-\varepsilon)(1+2\varepsilon)\big) = -\varepsilon + 2\varepsilon^2$ and $r_1 = \text{fma}\big(0 - \varepsilon(1+2\varepsilon)\big) = -\varepsilon - 2\varepsilon^2$. Thus $z_1 = c_0 = -2\varepsilon$, exactly $\Delta z_0 = z - z_0$.

Performance of twofold dividing is determined by two dividing operations, which are slow. This way, twofold "fast" dividing must be nearly 2x slower than standard, and so meet our 2x and 10x criteria. Testing shows this projection is right for Intel AVX processor, see <u>Performance</u> section below.

## Root

Let us compute square root of twofold argument. We assume approximate $\mathrm{sqrt}(x)$ function is available for a floating-point argument. Given $x = x_0 + x_1$, we let $z_0 = \mathrm{sqrt}(x_0)$, and we need $z_1$ approximating $\Delta z_0 = z - z_0$, where $z = \sqrt{x}$.

Let us start with simple case if $x_0 + x_1$ is non-overlapping "coupled". Here we can do Newton iterations. Because a Newton iteration nearly duplicates accuracy, one iteration is enough, assuming $z_0 = \mathrm{sqrt}(x_0)$ is accurate, maybe correctly rounded like IEEE-754-2008 standard requires.

If we utilize FMA for appropriate $x - z^2$, formula is:

$$z_1 \approx (x_0 + x_1 - z_0^2)/2z_0 \approx (x_1 + \mathrm{fma}(x_0 - z_0^2))/2z_0$$

Let us write this algorithm explicitly:

> <u>Algorithm RPF 1</u>. $\mathrm{sqrt}(x_0 + x_1) \to z_0 + z_1$
> (1) $z_0 = \mathrm{sqrt}(x_0)$
> (2) $z_1 = (x_1 + \mathrm{fma}(x_0 - z_0^2))/2z_0$

Note how this algorithm processes input below zero. Because $x_0 + x_1$ is non-overlapping, $x < 0$ implies $x_0 < 0$, and therefore $\mathrm{sqrt}(x_0)$ raises the domain error, which you can process later if necessary.

Simplifying for special case if $x_1 = 0$:

> <u>Algorithm RPF 1.1</u>. $\mathrm{sqrt}(x_0) \to z_0 + z_1$
> (1) $z_0 = \mathrm{sqrt}(x_0)$
> (2) $z_1 = \mathrm{fma}(x_0 - z_0^2)/2z_0$

Provided $\mathrm{sqrt}(x_0)$ is correctly rounded, RPF 1.1 must return non-overlapping "coupled" $z_0 + z_1$. Indeed, by design $z_1 \approx \Delta z_0$, while $|\Delta z_0| \leq \mathrm{ulp}(z_0)/2$ if library function's result is rounded correctly.

Now consider more complicated case if $x_0 + x_1$ is not "coupled". Here we cannot use Newton iterations that easily, because we cannot claim if $x_1$ is small comparing $x_0$. Our trick is reducing problem to known case of "coupled", renormalizing the input:

- $x_0 + x_1 \to u_0 + u_1$
- $\mathrm{sqrt}(u_0 + u_1) \to v_0 + v_1$
- $z_1 = (v_0 + v_1) - z_0$

For computing the difference $(v_0 + v_1) - z_0$ we can use twofold fast subtraction algorithm STF 1.1 and distil its result, which supplies very accurate result of subtraction. Explicitly, this algorithm is:

> <u>Algorithm RTF 1</u>. $\mathrm{sqrt}(x_0 + x_1) \to z_0 + z_1$
> (1) $z_0 = \mathrm{sqrt}(x_0)$
> (2) $x_0 + x_1 \to u_0 + u_1$
> (3) $\mathrm{sqrt}(u_0 + u_1) \to v_0 + v_1$
> (4) $(v_0 + v_1) \ominus z_0 \to w_0 + w_1$
> (5) $z_1 = w_0 \oplus w_1$

Special case if $x_1 = 0$ obviously reduces to RPF 1.1

Let us demonstrate how RTF 1 works if twofold $x_0 + x_1$ is not "coupled":

<u>Example</u>: Let $x_0 = x_1 = 1$. Here $z_0 = \mathrm{sqrt}(1) = 1$. Then $u_0 = 2$ and $u_1 = 0$, thus $v_0 + v_1 \approx \sqrt{2}$. Thus $z_1 \approx w_0 + w_1 \approx \sqrt{2} - 1$. Such $z_0 + z_1$ is not 2x-precise but still approximates $\sqrt{2}$ as expected.

Performance is determined by additional dividing and twofold add/subtract operations. If $\mathrm{sqrt}(x)$ itself is at least 10x slower than basic add/subtract, RPF/RTF must be 2-3 times slower than square root, thus must meet the "slower not more than 10x times if data in cache" performance criterion.

Concerning "2x slower for data in RAM" criterion, see our testing results in [Performance](#) section below.

## 2x

We cannot recommend twofold arithmetic for simulating 2x-precise calculations, it is not strict enough.

Exception is using 2x-precise twofold accumulator for summing a series of 1x-precise (regular) numbers. Adding a number to a twofold is strict, see comments to Algorithms ATS/STS 1.1

Anyway, if standard quad-precision type not available, or works too slowly, you may use twofolds over standard double as surrogate for quad. Worth using special function variants for "coupled" input, and renormalize output immediately. This would literally repeat the approach by Dekker [1].

For renormalizing immediately, you may use Fast-Two-Sum that costs only 3 add/subtract operations.

Renormalizing does not look to improve accuracy. However, consider sample long chain of calculations: adding unity to twofold accumulator, $2^{48}$ times basing on binary24 format:

- W/o renormalizing: result would saturate at $z_0 = z_1 = 2^{24}$, very far from correct
- With renormalizing: result would equal $z_0 = 2^{48}$ and $z_1 = 0$, exactly as expected

Note that you do not need to renormalize result in case if input $x_1 = y_1 = 0$, because output would be already non-overlapping "coupled". For add/subtract and multiply, this would be the exact transforms. For dividing and square root, see comments to these functions above.

# C++ test code

- Code archive
- Architecture

## Code archive

The performance test code, build/test scripts, and testing results are available free at our Web site [7]. There you can download and unpack code archive as zip-file: Twofold fast arithmetic, code.zip

This archive includes the folder named "code" with following sub-folders, each containing C++ sources, corresponding make file, and testing logs:

- applications
- perftest
- twofold

The folder "twofold" contains our experimental implementation of twofold "fast" arithmetic, and sanity test for it. Make file designed for Microsoft and GNU compilers; you can run it from command-line with make or nmake utility:

```
make gcc
nmake cl
```

The folder "perftest" contains the performance test, the universal make file for MS/GNU compiler, and testing results discussed in the next section named Performance.

Folder "applications" contains examples of using twofolds, discussed in section Applications below.

## Architecture

Here we briefly explain the coding style for better understanding the performance testing results. Main part is implementation of twofold algorithms, which directly encodes algorithms from [Background] and [Algorithms] sections above.

For best of SIMD performance, code utilizes AVX intrinsic provided with Intel, GNU, and Microsoft C++ compilers. To leverage C++ templates, we use very thin unified interface for basic arithmetic operations add, subtract, multiply, etc. for standard `float` or `double` and for AVX intrinsic types `__m256` or `__m256d`. This allows the same C++ code to target both vector (SIMD) and scalar data types.

<u>Fragment 1</u>: Uniform vector/scalar abstraction for fused-multiply-add

```cpp
#include <cmath> // scalar fma(x,y,z)
inline __m256d fma(__m256d x, __m256d y, __m256d z) { return _mm256_fmadd_pd(x,y,z); }
inline __m256  fma(__m256  x, __m256  y, __m256  z) { return _mm256_fmadd_ps(x,y,z); }
```

The following code fragment shows the twofold and "coupled" data types. Types are generic, assuming C/C++ standard `float` or `double` or AVX intrinsic type `__m256` or `__m256d` as a "number". Note the types hierarchy; we can assign a `coupled<T>` value to `twofold<T>` variable but not conversely:

<u>Fragment 2</u>: Twofold and couple-length "numbers"

```cpp
// Assume number is scalar single or double by IEEE-754,
// or vectored __m256 or __m256d of Intel AVX intrinsic:
template<typename number> struct twofold { number value, error; };
template<typename number> struct coupled: public twofold<number> {};
```

The arithmetic algorithms implemented as inline functions for best of compiler optimization. Note that strict-math compilation mode is required, as fast-math optimizations may eliminate the rounding tricks on which exact transforms base. The following fragment displays the two-product algorithm. Note, that we significantly use fused-multiply-add (FMA) here, so need a processor that supports fast FMA.

<u>Fragment 3</u>: Two-product algorithm (see Algorithm 3 from [Background])

```cpp
// Use fmadd(), so additional operation for negating:
template<typename T> inline coupled<T> pmul(T x, T y) {
        coupled<T> z;
        z.value = mul(x,y);
        z.error = fma(x,y,neg(z.value));
        return z;
}
```

The following code fragment implements twofold division if arguments are "dotted" numbers (not shaped as twofold or "coupled"). Note that the output is non-overlapped "coupled":

<u>Fragment 4</u>: Twofold division (Algorithms DTF/DPF 1.1.1)

```cpp
// Twofold divide, both x and y are dotted, so z is coupled:
template<typename T> inline coupled<T> tdiv(T x, T y) {
        T q0, q1, r1;
        q0 = div(x,y);          // q0 = x / y
        r1 = fma(neg(q0),y,x);  // r1 = x - q0*y
        q1 = div(r1,y);         // q1 = r1 / y
        coupled<T> z;
        z.value = q0;
        z.error = q1;
        return z;
}

// Coupled divide, x and y dotted:
```

```cpp
template<typename T> coupled<T> pdiv(T x, T y) {
        return tdiv(x,y);
}
```

The function names `tdiv`/`pdiv` overloaded for dotted and "shaped" arguments. Following is full list of functions for dividing. Lists for other functions look similarly:

Fragment 5: Twofold/coupled dividing interface

```cpp
template<typename T> twofold<T> tdiv(twofold<T> x, twofold<T> y);
template<typename T> twofold<T> tdiv(twofold<T> x,          T  y);
template<typename T> twofold<T> tdiv(          T  x, twofold<T> y);
template<typename T> twofold<T> tdiv(coupled<T> x, coupled<T> y);
template<typename T> twofold<T> tdiv(coupled<T> x,          T  y);
template<typename T> twofold<T> tdiv(          T  x, coupled<T> y);
template<typename T> coupled<T> tdiv(          T  x,          T  y);

template<typename T> coupled<T> pdiv(coupled<T> x, coupled<T> y);
template<typename T> coupled<T> pdiv(coupled<T> x,          T  y);
template<typename T> coupled<T> pdiv(          T  x, coupled<T> y);
template<typename T> coupled<T> pdiv(          T  x,          T  y);
```

The prefix "t" in the function name means twofold and prefix "p" means "coupled" type of output. Main set of functions implement algorithms for twofold arguments, plus the special algorithms in case if input is non-overlapping "coupled" or just a dotted number. Additional algorithms simulate Dekker arithmetic over "coupled" inputs by renormalizing the output, so ensuring result is also "coupled".

Generic type T may be scalar double/float or AVX vector __m256d/__m256. Even if T is vector, we add yet another vectoring level and define functions of array arguments. The tested compilers are very good in optimizing array functions, so we can utilize up to 90% percent of processor peak performance. Our test iterates the array calculations and measures the performance.

We add prefix "v" to the vector function names, and suffix "2" or "1" to distinguish functions with two or one twofold/coupled arguments. Resulting vector `r[]` is always shaped (twofold or coupled):

Fragment 6: Example of the array function interfaces:

```cpp
void vtadd2 (int m, twofold<__m256d> x[], twofold<__m256d> y[], twofold<__m256d> r[]);
void vtadd1 (int m, twofold<__m256d> x[],          __m256d  y[], twofold<__m256d> r[]);
void vtadd  (int m,          __m256d  x[],          __m256d  y[], twofold<__m256d> r[]);

void vpadd2 (int m, coupled<__m256d> x[], coupled<__m256d> y[], coupled<__m256d> r[]);
void vpadd1 (int m, coupled<__m256d> x[],          __m256d  y[], coupled<__m256d> r[]);
void vpadd  (int m,          __m256d  x[],          __m256d  y[], coupled<__m256d> r[]);
```

Functions tadd() of coupled arguments are tested only indirectly via calling from padd().

We measure performance in millions of twofold/coupled outputs per second. We call this metric mega-mega-operations-per-second, briefly mega-ops, or mega-flops. Our goal: be not slower than 1/10 of processor peak if data fit into CPU cache. If CPU peak were 10+ gigaflops for standard double type, our target is 1+ gigaflops for twofolds over doubles.

# Performance

## Test system

For performance testing we used laptop, built on ultra-low-voltage processor, with peak performance around 10 gigaflops for standard double type. Not very best choice for high-performance computing. But we still can use it for proving the concept, if twofold arithmetic at all can operate at 1+ gigaflops with this sort of modern processors, and that performance is slower by only 2x if data not in cache.

We have filtered intermittent effects in testing by repeating test runs and selecting lower-level results. This stabilize results and this way allows comparing performance of twofold arithmetic and of regular dotted floating-point operations.

The test system was HP Pavilion 15 laptop, built on Intel Core i5-4200U (Haswell) processor of nominal frequency 1.6 GHz and up to 2.6 GHz in turbo mode, the memory was 2x4 GB banks of PC-12800 (DDR3) so enabling up to 25.6 GB/s in overall. The compilers were GNU g++ 4.8.2 (Cygwin) and Microsoft Visual Studio 2013 Express. Performance data hereon are for GNU compiler that shows better megaflops.

The following table shows results of memory reading and copying and of dotted arithmetic. This result is for "vector" test with `__m256` and `__m256d` as basic types; see other data at my Web site [7]. CPU actual frequency was around 2.25 GHz as I could observe with Windows Task Manager. Note that CPU did a bit more, about 2.5 GHz in "scalar" test with `float` or `double` as basic types.

The function `vmem()` reads data from two arrays `x[]` and `y[]` and writes to another third array `r[]`, so simulating fetching data for arithmetic operations and storing results. The functions `vadd()`, `vmul()`, `vdiv()`, and `vsqrt()` actually perform the arithmetic operations:

Table 1: Memory copy and arithmetic over dotted data of "vector" basic types

| func | float | | | double | | |
|---|---|---|---|---|---|---|
| | small | medium | large | small | medium | large |
| vmem | 19753.5 | 2340.6 | 969.72 | 9924.99 | 1167.57 | 484.942 |
| vadd | 14774.8 | 2093.19 | 971.261 | 7398.14 | 1038.89 | 484.135 |
| vmul | 14644 | 2074.52 | 971.51 | 7228.24 | 1035 | 485.604 |
| vdiv | 1422.15 | 1409.71 | 999 | 363.132 | 361.956 | 356.078 |
| vsqrt | 1423.99 | 1421.31 | 1200.21 | 362.533 | 363.717 | 361.607 |

The test tries small, medium, and large arrays of around 100, 10 thousands, and 1 million of `float` or `double` elements, so the arrays fit the fastest L1 cache, fit the last-level (L3) cache, or do not fit CPU cache. The left part of the table is for single and the right is for double-precision basis type, or for the `__m256` and `__m256d` in "vector" test (which are the AVX 256-bit packs of 8 floats or of 4 doubles).

For add/subtract and multiply operations, this processor peak is around 18 gigaflops in single precision and 9 gigaflops in double if operating at 2.25 GHz, or 20 and 10 gigaflops if at 2.5 GHz. As we observed, CPU operated at nearly 2.25 GHz with twofold/coupled arithmetic, and at nearly 2.5 GHz with "scalar" dotted arithmetic testing. Such system's behavior looks caused by automatic balancing of CPU heating.

L1 cache performance looks enough to feed arithmetic if small arrays. Performance with large arrays limited by memory bandwidth, except very slow dividing square root of double-precision numbers.

With large arrays, performance looks 2x below our expectation. With ~1 gigaflops at single precision, each operation gets two numbers and writes one, so such performance implies reading 8 GB/s and writing 4 GB/s, so transferring 12 GB/s in overall. This is around half of bandwidth we expected.

If we look at results in more details available at our Web site [7], there we can see twice-higher results intermittently occurring in test runs. We filtered such full-bandwidth results away, and analyzed easier to reproduce half-bandwidth results. Note that achieving maximal CPU performance is not a subject of this work, enough if we can compare results for dotted and twofold/coupled arithmetic.

## Test results

Look at performance results for twofold arithmetic over the "vector" __m256 and __m256d types:

<u>Table 2</u>: Twofold performance over AVX "vector" types

| func | float | | | double | | |
|---|---|---|---|---|---|---|
| | small | medium | large | small | medium | large |
| vtadd2 | 2238.72 | 1021.61 | 459.423 | 1141.65 | 509.554 | 230.971 |
| vtadd1 | 2495.9 | 1142.39 | 509.304 | 1265.07 | 570.755 | 259.465 |
| vtadd | 2926.55 | 1337.6 | 618.993 | 1456.83 | 669.919 | 306.501 |
| vtsub2 | 2183.93 | 976.855 | 429.351 | 1122.57 | 495.039 | 223.265 |
| vtsub1 | 2518.56 | 1136.72 | 514.479 | 1268.57 | 565.852 | 263.608 |
| vtsub | 2973.46 | 1331.08 | 661.223 | 1510.45 | 691.52 | 329.568 |
| vtmul2 | 3653.93 | 1036.89 | 469.32 | 1830.35 | 518.363 | 237.62 |
| vtmul1 | 4891.44 | 1193.9 | 553.781 | 2445.97 | 596.938 | 279.025 |
| vtmul | 6115.15 | 1398.88 | 664.782 | 2774.36 | 699.494 | 331.331 |
| vtdiv2 | 642.519 | 641.706 | 483.259 | 162.599 | 162.365 | 160.154 |
| vtdiv1 | 637.749 | 636.9 | 564.571 | 162.347 | 162.279 | 160.605 |
| vtdiv | 637.58 | 638.315 | 590.21 | 162.37 | 162.256 | 158.014 |
| vtsqrt1 | 430.109 | 429.065 | 382.773 | 108.182 | 107.982 | 106.058 |
| vtsqrt | 635.026 | 625.015 | 617.678 | 162.262 | 161.748 | 160.645 |

Twofold add/subtract functions operate at 15% of actual peak performance for dotted as measured by the dot-test (Table 1), show 1120+ versus 7400 megaflops for double precision if small data in L1 cache. If one or both arguments are dotted, performance is 17% and 20% correspondingly. Versus theoretical peak, which is 9 gigaflops for doubles at 2.25 GHz, twofold add/subtract performs at 12-17% of peak.

Twofold multiply operates at 24% to 38% versus measured actual peak of 7300 megaflops for doubles. Versus theoretical peak of 9 gigaflops at 2.25 GHz, twofold multiplication operates at 20% to 30%.

For large data not in cache, twofold add/subtract/multiply performance looks driven by memory, and appears nearly 2.1 times slower than dotted, 230 twofold versus 484 dotted megaflops for doubles if both arguments are shaped.

In overall, in terms of our 10x and 2x criteria:

- twofold add/subtract/multiply meet the "10x" criterion, even with significant handicap
- nearly meet the "2x" criterion, though not quite, twofold is 2.1 times slower than dotted

Twofold dividing and square root appear 2.2x and 3.3x slower than dotted, as we presumed.

For more data, please see our Web site [7]. There you can find results for "coupled" functions, which appear nearly same as for twofold. And results for scalar types with compiler-driven vectoring for AVX. The scalar results look fine for dotted functions, but look really awful for twofold/coupled. Thus, we need to optimize manually for good performance with twofold/coupled arithmetic.

## Conclusion

Twofold arithmetic can operate at 10% to 40% of processor peak if data in cache. Performance with large data is nearly 2x times below the dotted arithmetic, as limited with memory bandwidth.

This way, this implementation of twofold arithmetic meets our 10x performance criterion with good handicap, and almost meets the 2x criterion. This is very good result: we meet our performance goal!

We need to utilize 100% of processor's peak performance to meet the performance goals with twofolds. Manual code vectoring is necessary here; automatic vectoring by modern C++ compilers is not enough.

## Applications

- [Summation](#)
- [Linear Ax=f](#)
- [Square root](#)
- [Conclusion](#)

## Summation

Simplest but important test is just summation $s = \sum x_i$ or scalar product $p = \sum x_i y_i$. Let us consider even simpler partial case (I heard it from Marius Cornea). Imagine a counter $s$ designed to accumulate sum of many identical values $x_i = t$, where both $s$ and $t$ are not precise. This example is from real tech accident, where $s$ was counting time (in seconds) and $t$ was equal to 1/10 second.

In that accidental case, $s$ and $t$ were implemented as single-precision floating-point values, which allows relative inaccuracy $\varepsilon \sim 10^{-7}$ that seems good enough. The problem was inaccuracy accumulation in $s$ if a lot of summations. For example, if we wait 100 hours, the accumulated error in $s$ would get as large as almost 4 hours or 4%, and the error would exceed 40% of correct value if we wait 1000 hours.

Let us show how twofolds could address this situation. Following is testing log that you can find under "applications" folder in our code archive (see [C++ test code](#) section above). The test named `test100h` shows the floating-point values of $s$ and $t$ if implemented as twofold over `float` or `double` type. Here the twofold value of "1/10 s" refers to $t$, and "result" corresponds to $s$:

```
test: type=float, hours=100
    1/10 s: 0.1[-1.49012e-09]
    result: 96.3958[3.54008] hours
    expect: 100 hours
test: type=double, hours=100
    1/10 s: 0.1[0]
    result: 100[3.33695e-09] hours
    expect: 100 hours
test: type=float, hours=1000
    1/10 s: 0.1[-1.49012e-09]
    result: 582.542[461.249] hours
    expect: 1000 hours
test: type=double, hours=1000
    1/10 s: 0.1[0]
    result: 1000[-6.12184e-07] hours
    expect: 1000 hours
```

Here each number in square brackets is the error part of the corresponding twofold number.

No surprise, double precision appears enough to keep accumulated error reasonably low, and twofold arithmetic shows specific estimate: around $3_{10}-9$ for easier 100-hours case, and $6_{10}-7$ for 1000 hours.

For single precision, twofold arithmetic provides much worse estimates for accumulated errors: around 3.5 for 100-hours test, and around 461.2 per 1000 hours. Such estimates correctly signal on the problem with awful main results: 96.4 instead of 100 hours, and 582.5 instead of 1000. With such estimates, the controlling software could stop operations and so prevent the accident caused by the summation error.

## Linear Ax=f

Consider using twofolds in direct Gauss solver for $Ax = f$ linear system with very simple matrix $A$, which would be just a Jordan block with a small $\lambda$:

$$A = \begin{bmatrix} \lambda & 1 & 0 \\ 0 & \ddots & 1 \\ 0 & 0 & \lambda \end{bmatrix}$$

Well known that such simple system might be very bad for solving numerically if $\lambda$ is far from unity. Let us see if twofolds can identify and assess the accumulation of inaccuracy in such numeric solution.

Following is testing log for the "gauss" test you can find under the "applications" folder in the twofolds code archive (see section C++ test code above). Here we solve with 3x3 matrix, and $\lambda$ equal to 1/10 for well-conditioned case and 1/1000 for ill-conditioned. The system's right part $f$ especially designed for simple expected solution $x$.

The test solves the system with direct Gauss method. Test calculates with twofolds over `float` or `double` type. For initializing matrix $A$, we convert $\lambda$ to `double` type and then to `float`, so diagonal elements for `float` case include rounding error estimates.

Following piece of test log is for well-conditioned 3x3 system:

```
test, float, well3
    A
        0.1[-1.49012e-09]    1[0]              0[0]
        0[0]                 0.1[-1.49012e-09] 1[0]
        0[0]                 0[0]              0.1[-1.49012e-09]
    f
        11[0]    11[0]    1[0]
    x (expected)
        10[0]    10[0]    10[0]
    x (solution)
        10[0]    10[0]    10[0]

test, double, well3
    A
        0.1[0]    1[0]     0[0]
        0[0]      0.1[0]   1[0]
        0[0]      0[0]     0.1[0]
    f
        11[0]    11[0]    1[0]
    x (expected)
        10[0]    10[0]    10[0]
    x (solution)
        10[-5.05151e-14]    10[4.996e-15]    10[-5.55112e-16]
```

Occasionally, twofolds fail identifying any error in `float` calculations. But twofold estimates look well found for `double` case: twofold says, error of x[2] is around $5_{10}\text{-}16$, around 10 of ULP, like we should expect. And error gradually increases, by $1/\lambda = 10$ times in x[1], and then by 10 times more in x[0].

Following piece of test log is for ill-conditioned 3x3 matrix:

```
test, float, ill3
    A
        0.001[-4.74975e-11]    1[0]                0[0]
        0[0]                   0.001[-4.74975e-11] 1[0]
        0[0]                   0[0]                0.001[-4.74975e-11]
    f
        1001[0] 1001[0] 1[0]
    x (expected)
        1000[0] 1000[0] 1000[0]
```

```
        x (solution)
            939.026[60.9742]     1000.06[-0.0609741]     1000[6.10351e-05]

test, double, ill3
    A
        0.001[0]      1[0]        0[0]
        0[0]      0.001[0]        1[0]
        0[0]          0[0]    0.001[0]
    f
        1001[0] 1001[0] 1[0]
    x (expected)
        1000[0] 1000[0] 1000[0]
    x (solution)
        1000[-2.07959e-08]     1000[2.07959e-11]     1000[-2.08167e-14]
```

Here, inaccuracy gradually increases by $1/\lambda = 1000$ times from x[2] to x[1] and then to x[0]. Twofolds correctly assess this accumulated inaccuracy for both `float` and `double` types. For `float` case, solution looks obviously wrong in x[0], and twofold correctly shows how much wrong is it.

## Square root

To illustrate twofold square root, we solve quadratic equation $ax^2 + bx + c = 0$ with school formula:

$$x_{0,1} = \frac{-b \pm \sqrt{b^2 - 4ac}}{2a}$$

We let $a = 1$ and $b = 2$, and try $c$ very close to 0 or to 1, specifically $c = 10^{-8}$ or $c = 1 \pm 10^{-8}$. First case examines accuracy loss in $-b + \sqrt{b^2 - 4ac}$, and second is about square root of inaccurate input.

Note that binary32 type (aka `float`) cannot represent $1 \pm 10^{-8}$, thus $b^2 - 4ac$ would be wrong with single precision. For `double` type, this formula must cause losing around half of significant digits.

Following is fragment from test log found at our Web site [7]:

```
test: type=float
    a: 1[0]
    b: 2[0]
    c: 1e-08[6.07747e-17]
    d: 2[1e-08]
    x0: -2[-5e-09]
    x1: 0[5e-09]
test: type=double
    a: 1[0]
    b: 2[0]
    c: 1e-08[0]
    d: 2[6.07747e-17]
    x0: -2[-1.4141e-16]
    x1: -5e-09[3.03874e-17]
```

Here we let $c = 10^{-8}$, so discriminant $d = \sqrt{b^2 - 4ac}$ must equal $\sqrt{3.99999996} \approx 2 - 10^{-8}$, and result must be $x_0 \approx -1.999999995$ and $x_1 \approx -5.0000000125 \cdot 10^{-9}$.

Double precision results fit these expectations modulo printing fewer significant digits; and accuracy estimate for $x_1$ predictably says around half of significant digits is lost. Single precision (`float`) is not enough, and twofold correctly assesses inaccuracy; particularly, $x_1$ looks completely wrong here.

Note how twofold square root behaves if argument gets out of range. Here we let $c = 1 + 10^{-8}$, so result is NaN, which we can identify or propagate through further computations:

```
test: type=float
    a: 1[0]
    b: 2[0]
    c: 1[1e-08]
    d: 0[nan]
    x0: -1[nan]
    x1: -1[nan]
test: type=double
    a: 1[0]
    b: 2[0]
    c: 1[0]
    d: nan[nan]
    x0: nan[nan]
    x1: nan[nan]
```

## Conclusion

These are very simple though typical examples of how twofolds could assure more reliable computing.

Note that in these tests, twofolds would not over-estimate actual inaccuracy; avoid paranoid signaling on problem if situation is actually fine. Such good behavior is not just occasion:

If basic precision is enough, like for `double` here, 2x-higher precision of twofold is moreover enough, and $z_1$ is relatively small. But if $z_1$ is large, this almost for sure means basic precision is not enough.

This is somewhat similar to memory parity check: if parity bit looks fine this means nothing as system might just miss a problem, but if parity bit is wrong this must signal on a real problem.

Thus, potentially you may double-check almost any math result with minimal risk of irrelevant panic.

## Discussion

- [Motivation](#)
- [Technology](#)
- [Processors](#)
- [Compilers](#)

### Motivation

Our motivation follows simple philosophy that isolating arithmetic difficulties is job for technology. Human must have privilege be unaware on too many details, and concentrate on areas of interest, construction, science, education, etc.

We design a twofold daemon to check accuracy of all floating-point results. Daemon would work in on-fly manner, in parallel with main computations. Daemon would avoid needless panic, and its cost would be affordable for majority of mathematic computations.

Checking all results with twofolds would simplify programming mathematics, let human encode math formulas directly "as is", and fix only in case of problems. If no problems found, such checking anyway assures reliability of math computations.

### Technology

For implementing the twofolds daemon, we adapt Dekker arithmetic [1] for modern processors, which were not available in 1970[th]. This allows twofolds be very fast, so the daemon would minimally damage overall performance of your program.

Computers evolve very quickly duplicating capacity every 18 months according to Moore's law. So cost of checking would seem negligible very soon, while benefits are substantial. Besides, future computers could support twofolds in hardware, so minimizing burden.

Application examples show that twofolds do not tend overestimating problems. This good property is not occasional but follows from twofold's approach. We use 2x-higher precision for assessing accuracy of 1x-precise results. If 1x-precision is fine, then 2x is fine moreover and assessment converges.

This potentially allows checking all sorts of mathematic calculations with minimal risk of irrelevant panic. If we instrument any mathematic program by replacing all `float` and `double` variables with twofolds, such program would compute bitwise same result plus check its accuracy.

### Processors

How difficult might be supporting twofold and similar techniques (like "double-double") in future CPUs?

We need faster Two-Sum and Two-Product operations described in [Background](#). Two-sum can be as fast as three add/subtract operations according to Algorithm 1.1, so results could retire every 3 CPU ticks. If a more expensive implementation with a conveyor, ultimately results could retire every 1 tick.

Thus, twofold summation Algorithms ATF/STF 1 could take only 5 or even 3 ticks instead 8, doing around twice faster so far. Special summation with Algorithm ATF 1.1, could cost 4 or even 2 ticks instead of 7, so operate 2x or even 3x faster in important cases of sum $\sum x_i$ and dot-product $\sum x_i y_i$.

Two-Product is easier than FMA, so can be same fast, retire every tick. With such improvement, twofold multiplication Algorithm MTF 1 could cost 7 ticks instead of 8, which speedup does not look critical.

### Compilers

For best performance, twofold/coupled arithmetic needs support in compilers. Twofold arithmetic formulas seem too complicated for automatic vectoring, so compilers should learn these patterns.

## Gratitude

Here we would like to thank the following people with whom we have discussed these ideas and results:

- Marius Cornea (Intel)
- Bob Hanek (Intel)
- Victor Kostin (Intel)
- Dmitry Baksheev (Intel)
- Evgeny Petrov (Intel)
- Alexander Semenov (UniPro)
- Ivan Golosov (UniPro)

Here, "Intel" is Intel Corporation (http://intel.com), UniPro web site is http://unipro.ru/eng/index.html

## References

[1] T. Dekker, A Floating-Point Technique for Extending the Available Precision, Numer. Math. 18, 224-242 (1971)

[2] D. Priest, On Properties of Floating Point Arithmetics: Numerical Stability and the Cost of Accurate Computations, ftp://ftp.icsi.berkeley.edu/pub/theory/priest-thesis.ps.Z

[3] J. Shewchuk. Adaptive precision floating-point arithmetic and fast robust geometric predicates. Discrete & Computational Geometry, 305–363, 1997.

[4] Y. Hida, X. Li, D. Bailey, Library for Double-Double and Quad-Double Arithmetic, http://web.mit.edu/tabbott/Public/quaddouble-debian/qd-2.3.4-old/docs/qd.pdf

[5] High-Precision Software Directory by D. Bailey et al, http://crd-legacy.lbl.gov/~dhbailey/mpdist/

[6] J.-M. Muller, et al. Handbook of Floating-Point Arithmetic. Springer, 2010

[7] See more materials at our Web site: https://sites.google.com/site/yevgenylatkin/